\magnification=1200

\input amstex
\documentstyle{amsppt}

\NoBlackBoxes 
\NoRunningHeads
\nologo

\pagewidth{6.5truein}
\pageheight{9truein}


\topmatter
\title
Commutative Algebras of Ordinary 
Differential Operators with Matrix Coefficients
\endtitle
\author
Masato Kimura$^{1)}$ and Motohico Mulase$^{2)}$  \\
\endauthor
\affil
$^{1)}$Department of Mathematics\\
University of Wisconsin\\
Eau Claire, WI 54702\\
kimuram\@uwec.edu\\
\\
$^{2)}$Department of Mathematics\\
University of California\\
Davis, CA 95616\\
mulase\@math.ucdavis.edu\\
\endaffil
\abstract
A classification of commutative
integral domains consisting of ordinary 
differential operators with matrix coefficients 
is established in terms of morphisms 
between algebraic curves.
\endabstract
\thanks
$^{2)}$Research supported in part by 
NSF Grant DMS--94--04111.
\endthanks
\endtopmatter
\document


\def\C{{\Bbb C}}
\def\cfield{{\Bbb C((z))}}
\def\cring{{\Bbb C[[z]]}}
\def\fullloop{{gl\big(n,\Bbb C((z))\big)}}
\def\posloop{{gl\big(n,\Bbb C[[z]]\big)}}

\def\isom{\cong}
\def\tensor{\otimes}
\def\dsum{\oplus}
\def\rarrow{\rightarrow}

\def\longrarrow{\longrightarrow}

\def\Spec{{\text{\rm Spec}}}
\def\Proj{{\text{\rm Proj}}}
\def\End{{\text{\rm End}}}
\def\Ker{{\text{\rm Ker}}}

\def\Coker{{\text{\rm Coker}}}
\def\ord{{\text{\rm ord}}}
\def\rank{{\text{\rm rank}}}
\def\trace{{\text{\rm trace}}}
 
\def\hookrarrow{\hookrightarrow}

\bigskip\bigskip
\subhead
1. Introduction and the main results
\endsubhead
\medskip

The purpose of this paper is to establish a geometric
classification of  commutative integral domains
consisting of linear ordinary differential operators 
with matrix coefficients. {From} an integral domain
of matrix ordinary differential operators we construct
a morphism between two irreducible algebraic curves
and certain torsion-free sheaves on them 
which are compatible with the
morphism. Conversely, from a set of geometric data 
consisting of 
an arbitrary morphism of integral curves and torsion-free
sheaves on them with vanishing cohomology groups, we 
construct an integral domain  of matrix
ordinary differential operators. These two constructions
are inverse to one another. In this correspondence, arbitrary
integral curves (singular curves as well) appear.

It was G.~Wallenberg  who first recognized the
rich mathematical structure 
 in  commuting
ordinary differential operators. Through an
explicit computation he observed that the 
coefficients of two commuting linear ordinary differential
operators of orders 2 and 3 with scalar function
coefficients are given by elliptic
functions. The 1903 paper \cite{21} thus
initiated the long history of  attempts  towards
the classification of commutative subalgebras of 
the ring of ordinary differential operators with scalar
coefficients. The classification problem was finally
completed in 1990 by Mulase \cite{13} in terms of
the moduli spaces of pointed algebraic curves and vector
bundles defined on them. His results are based on the earlier work
of Schur \cite{18}, Burchnall and Chaundy 
\cite{3},
Gel'fand and Dickey \cite{5}, 
Krichever \cite{8},
Mumford \cite{15}, Verdier \cite{20}, 
Sato \cite{17} and others. Previato and Wilson 
\cite{16} also made a decisive contribution to this 
problem. 

In 1993, Adams and Bergvelt \cite{1}
discovered a beautiful three-fold
relation between maximally commutative subalgebras
of the loop algebra of $gl(n,\C)$, certain
morphisms between algebraic curves, 
and infinite-dimensional
integrable systems generalizing the KP equations. 
Their work gave a big impact to the
 further development of  Li and Mulase 
\cite{9} in establishing
 a characterization of the generalized 
Prym varieties in terms of the 
multi-component KP equations. As a byproduct,  a
new class of commutative algebras of ordinary
differential operators with matrix coefficients was
constructed from morphisms between algebraic 
curves.

In his dissertation research, Kimura \cite{7} realized 
that a more general class of commutative subalgebras
of the loop algebra than those
studied by \cite{1} and \cite{9}
would lead to morphisms
between singular algebraic curves. A similar
observation was made by Donagi and Markman
\cite{4}. 

It suggests the possibility of constructing a type of
{\it Galois theory\/} of commutative algebras of differential
operators. Let $B$ be a commutative 
algebra of ordinary differential
operators. If we choose an element $P\in B$, then 
$B$ is an extension of a subalgebra $\C[P]$. If 
these algebras
correspond to algebraic curves, then the inclusion relation
should correspond to a morphism between the curves.
Do we obtain an arbitrary morphism of curves in this way?
Conversely, if we start with an arbitrary morphism of
curves, can we construct a commutative algebra of 
ordinary differential operators and its extension 
that
correspond to the original geometric situation?

We answer these question in this paper. 
The complexity of nilpotent elements and
zero divisors arises in the theory of
differential operators with matrix coefficients. 
{From} a geometric point of view,
 it is natural to consider morphisms between 
irreducible curves as building blocks for more general 
cases. So we restrict ourselves to  irreducible curves,
and deal only with integral domains of differential
operators here. The  general case will be
treated elsewhere. 

We establish he following two theorems in this paper.

\medskip
\proclaim{Theorem 1.1}
Let $f:C\longrarrow C_d$ be an arbitrary morphism
of degree $n$
between reduced irreducible algebraic curves $C$ and
 $C_d$
defined over the field $\C$ of 
complex numbers, and let $p\in C_d$ be
a smooth  point
of $C_d$. Then for every torsion-free
rank one sheaf $\Cal L$ defined on $C$ such that
$$
H^0(C,\Cal L) = H^1(C,\Cal L) =0\;,
$$
there is an injective algebra homomorphism
$$
j_{\Cal L}:
H^0(C\setminus f^{-1}(p),\Cal O_C)\hookrarrow
gl(n,D)
$$
of the coordinate ring of $C$ into the ring $gl(n,D)$ of 
matrix ordinary differential operators,
where 
$$
D = \big(\C[[x]]\big)\left[\frac{d}{dx}\right]
$$
is the algebra of differential operators with 
coefficients in the ring $\C[[x]]$ of  formal power series 
in one variable.
\endproclaim
\medskip\noindent
The next theorem gives the converse.
\medskip
\proclaim{Theorem 1.2}
Let $B\subset gl(n,D)$ be a commutative 
integral domain consisting of 
matrix differential operators. We assume that 
 $B$ is elliptic in the sense of Verdier, namely, 
$B$ has a monic element of positive order, say
$$
P = I_n\cdot \left(\frac{d}{dx}\right)^r + {\text{  lower order
terms  }} \in B\;.
$$
Consider  a $\C$-subalgebra 
$$
B_d  \subset B\cap \C((P^{-1/r})) 
$$
of $B$.
We further impose that
\roster
\item The subalgebra $B_d$ has rank one, i.e., 
the greatest common divisor of the orders of operators
in $B_d$ is equal to $1$;
\item $B$ is a rank $n$ module over $B_d$. (Since $B$
is integral and hence
 torsion-free over $B_d$, we can define its rank.)
\endroster
Then there exist a degree $n$ morphism of integral algebraic curves
$$
f:C\longrarrow C_d\;,
$$
a torsion-free rank one sheaf $\Cal L$ defined over $C$, and
a smooth point $p\in C_d$ such that the construction
of Theorem~1.1 gives back  to $B$. Moreover,
the collection of the data $(B_d, B, P)$ satisfying the
above conditions and the set of geometric data
$\big\langle
f:C\rarrow C_d,\Cal L, p, z, \psi
\big\rangle$ satisfying the
conditions of Theorem~1.1 
are in one-to-one correspondence.

\endproclaim
\medskip
\demo{Remark 1}
The ground field should be an algebraically closed field of
characteristic 0. We take it to be just the field of complex 
numbers.
\enddemo
\demo{Remark 2}
If we start with a non-integral curve $C$ and a morphism onto 
an integral curve $C_d$ with the same conditions as in 
Theorem~1.1, then we obtain a non-integral commutative
subalgebra of $gl(n,D)$. 
\enddemo
\medskip
\noindent

This paper is organized as follows. In Section~2, we
study commutative subalgebras of a formal loop algebra
and give a method of constructing morphisms between
algebraic curves from the data of commutative loop
algebras. In Section~3, we establish the converse direction, 
i.e., from a set of geometric data consisting of a morphism
between integral curves and vector bundles on them, we 
give a set of algebraic data which are commutative subalgebras
of a loop algebra and a point of an infinite Grassmannian.
The relation between matrix ordinary differential operators
and the loop algebras is studied in Section~4. More detailed
study has been given in \cite{9}. In Section~5 we complete the
proof of Theorem~1.1.  Our theory includes the {\it 
spectral curves} of \cite{2}, \cite{7}, \cite{11} as
special cases. In Section~6, we show that from
each spectral curve we can construct a commutative
integral domain of matrix ordinary differential operators. 
The proof of Theorem~1.2 is completed in Section~7.
\bigskip\bigskip
\subhead
2. Formal loop algebras and algebraic curves
\endsubhead
\medskip

Let $\Bbb C$ be the field of complex numbers and
$\Bbb C[[z]]$ the ring of formal power series in one
variable $z$. The field of formal Laurent series,
$\Bbb C((z))$, is the field of fractions of $\Bbb C[[z]]$.
There is a $\C$-vector space direct sum decomposition
$$
\C((z)) = \C[z^{-1}] \dsum z\cdot \C[[z]]\;.
$$
The $z$-adic topology of $\cfield$ is introduced by
defining $z^k\cdot \cring$, $k\in \Bbb Z$, as the
system of open neighborhoods of $0\in \cfield$.
An element $a\in \cfield$ is said to have {\it order} $m$,
which is denoted by $\ord(a) = m$,
if $a\in z^{-m}\cdot \cring \setminus z^{-m+1}\cdot \cring$.
The negative sign indicates that we are considering the
{\it pole order\/} of elements of $\cfield$.

We fix a positive integer $n$ throughout the paper. Let us 
consider the formal loop algebra $gl\big(n,\cfield\big)$.
There is a natural filtration 
$$
\fullloop = \bigcup_{k\in \Bbb Z} z^k\cdot \posloop
$$
compatible with the addition and multiplication:
$$
\cases
+: z^k\cdot \posloop \times z^{k}\cdot \posloop
\longrarrow z^{k}\cdot \posloop\\
\cdot: z^k\cdot \posloop \times z^{\ell}\cdot \posloop
\longrarrow z^{k+\ell}\cdot \posloop\;.
\endcases
$$
The main object of this paper is a (commutative) integral
domain $A\subset \fullloop$. We denote
$$
A^{(k)} = A\cap \bigg(z^{-k}\cdot \posloop\bigg)\;.\tag 2.1
$$
Thus there is a natural filtration 
$$
A = \bigcup_{k\in\Bbb Z} A^{(k)}
$$
compatible with the addition and multiplication in $A$.
Let us choose a scalar diagonal $\C$-subalgebra $A_d$ of $A$:
$$
\C\subset A_d \subset A\cap \bigg(\cfield\cdot I_n\bigg)\;.
$$
We identify $A_d$  as a subalgebra of $\cfield$.  Since
$A_d\subset A$ and $A$ is integral, $A$ is a torsion-free
module over $A_d$. We define
$$
\rank(A_d) = GCD\{\ord(a)\;|\; a\in A_d\}\;.
$$
In order to avoid the extra complexity,
we assume that the commutative integral domain $A$ 
and a scalar diagonal subalgebra $A_d$ satisfy
the following conditions:
\medskip
\proclaim{Condition 2.1}
\roster
\item $\rank(A_d) = 1$.
\item As a torsion-free $A_d$-module, $A$ is of rank $n$.
\endroster
\endproclaim
\medskip
\proclaim{Theorem 2.2}
Let $(A_d, A)$ satisfy Condition~2.1. 
If $A$ satisfies further
that
$$
A^{(-1)} = A\cap \bigg(z\cdot
\posloop\bigg) =0\;,
$$
then the pair $(A_d,A)$ defines a morphism 
$$
f: C\longrarrow C_d
$$
of degree $n$ between two 
reduced irreducible complete algebraic curves
$C$ and $C_d$.
\endproclaim
\demo{Proof}
Let $A_d ^{(m)} = A_d \cap A^{(m)}$. We note that 
$$
\dim_{\C} A_d ^{(m)}\big/ A_d ^{(m-1)} \le 1\;. \tag 2.2
$$
Since $A_d$ is of rank 1, there are elements $a$ and $b$ in $A_d$ 
 such that $GCD(\ord(a),\ord(b)) = 1$. Let 
$$
N_{A_d} = \{\ord(h)\;|\; h\in A_d\} \subset \Bbb N = \{0,1,2,3,\cdots\}\;.
$$
Then it is easy to show that $N_{A_d}$ contains all integers
greater than
$$
\ord(a)\cdot \ord(b) - \ord(a) - \ord(b)\;.
$$
This fact, together with (2.2) and
 that $N_{A_d}$ does not contain any negative integers,
imply that 
\roster
\item $A_d$ is a finite module over $\C[a,b]$, and
\item $b$ is algebraic over $\C[a]$.
\endroster
Thus $\Spec(A_d)$ is an affine algebraic curve.
Following \cite{9}, \cite{12}, \cite{13} and \cite{15}, 
we define a graded algebra
$$
gr(A_d) = \bigoplus_{m=0} ^{\infty} A_d ^{(m)}
$$
and a projective scheme 
$$
C_d = \Proj(gr(A_d))\;.\tag 2.3
$$
Let
$$
A_d\left[\frac{1}{a}\right]_0
=\left. \left\{\frac{h}{a^m}\;\right|\;
m\ge 0, \; h\in A_d,\; \ord(h)\le m\cdot \ord(a)\right\}
\subset \cring\;.
$$
Then 
$$
A_d\cap A_d\left[\frac{1}{a}\right]_0 = A_d ^{(0)} = \C\;,
$$
because $A_d ^{(-1)} = 0$. Therefore, 
$$
\Proj(gr(A_d)) = \Spec(A_d) \cup 
\Spec\left(A_d\left[\frac{1}{a}\right]_0\right)\;, \tag 2.4
$$
because a regular function on $\Spec(A_d)$ that is also 
regular on $\Spec\left(A_d\left[\frac{1}{a}\right]_0\right)$
is a constant. Thus these two affine subschemes cover the
whole projective scheme $\Proj(gr(A_d))$. Since $\Spec(A_d)$
is an affine curve, $C_d$ is a complete algebraic curve.
As is shown in the literature cited above, $C_d$ is a
one-point completion of $\Spec(A_d)$. The $z$-adic
formal completion of 
$$
A_d\left[\frac{1}{a}\right]_0
\subset \cring
$$ 
is equal to $\cring$. The fact that
$$
A_d ^{(0)} = A_d \cap \cring = \C
$$
shows that 
$$
C_d = \Spec(A_d) \cup \{p\}\;, \tag 2.5
$$
where $p\in C_d$  is the unique nonsingular 
geometric point of $\Spec(\cring)$. Since it
is the point of $C_d$ defined by
an equation  $z=0$, $p$ is a (rational) nonsingular 
point of $C_d$. Thus $C_d$ is a one-point completion 
of $\Spec(A_d)$ by
a nonsingular point.

Similarly, we define
$$
gr(A) = \bigoplus_{m=0} ^{\infty} A^{(m)}\;.
$$
Since $gr(A)$ is a finite-rank module over $gr(A_d)$, 
$$
C = \Proj(gr(A)) \tag 2.6
$$
is a complete algebraic curve, and the natural inclusion
$$
gr(A_d) \subset gr(A)
$$
defines a surjective morphism
$$
f: C\longrarrow C_d\;. \tag 2.7
$$
Since $A$ is of rank $n$ over $A_d$, the degree of the morphism
$f$ is also $n$. This completes the proof.
\enddemo
\medskip
Note that $A^{(0)}$ is a $\C$-subalgebra of $A$. It is an integral
domain, and since $A^{(-1)}=0$, it is a subalgebra of
$gl(n,\C)$. Therefore, every element of $A^{(0)}$ is algebraic
over $\C$, and hence
$$
A^{(0)} = \C\;.
$$
Let us define
$$
A\left[\frac{1}{a}\right]_0
= \left. \left\{\frac{h}{a^m}\;\right|\;m
\ge 0, \; h\in A^{(m\cdot \ord(a))}\right\}
\subset \posloop\;.
$$
Then
$$
A \cap A\left[\frac{1}{a}\right]_0 = A^{(0)} = \C\;.
$$
Thus
$$
C = \Proj(gr(A)) = \Spec(A)\cup \Spec\left(
A\left[\frac{1}{a}\right]_0\right)
$$
because of the same reason as in (2.4). However, this time
$C$ is neither a one-point completion of $\Spec(A)$ nor 
a completion by smooth points.
The $z$-adic formal completion of $A\left[\frac{1}{a}\right]_0$
is 
$$
A\left[\frac{1}{a}\right]_0 ^{\wedge} = 
A\left[\frac{1}{a}\right]_0 
\widehat{\bigotimes}_{A_d\left[\frac{1}{a}\right]_0} \cring\;,
$$
which is a module of rank $n$ over
$\cring$ because of the second item of Condition~2.1.
The formal scheme 
$$
\Spec\left(A\left[\frac{1}{a}\right]_0
^{\wedge}\right) = \widehat{C}_{f^{-1}(p)}
$$
is the formal completion 
of $C$ along the divisor
$f^{-1}(p)$. In general, it is not smooth.
If we consider even a larger extension
$$
A\left[\frac{1}{a}\right]_0 ^{\wedge}
{\bigotimes}_{\cring} \cfield\;,
$$
then it is an $n$-dimensional vector space over the
field $\cfield$. We can choose a $\cfield$-linear basis
for the above vector space from $A\left[\frac{1}{a}\right]_0$,
say $X_1, \cdots,  X_n$.
Then 
$$
\cring\big[X_1,\cdots,X_n\big] \subset 
A\left[\frac{1}{a}\right]_0 ^{\wedge}\;.
$$
If there is an element $X\in A\left[\frac{1}{a}\right]_0$ such that
$$
\cring [X] = A\left[\frac{1}{a}\right]_0 ^{\wedge}\;,
$$
then
$$
\widehat{C}_{f^{-1}(p)} = \Spec(\cring [X])
= \Spec\left(\frac{\cring[t]}{(ch_X(t))}\right)\;,
$$
where $ch_X(t)\in\cring [t]$
is the characteristic polynomial of $X$.
This is the case of a {\it spectral curve\/}, where the
covering of $\Spec(\cring)$ at infinity is defined by the
characteristic polynomial $ch_X(t)$. It can well be 
singular, but the singularity of $\widehat{C}_{f^{-1}(p)}$
is far more general than those appearing in the case of
spectral curves.

\bigskip\bigskip
\subhead
3. The geometric data and the Grassmannian
\endsubhead
\medskip

The geometric data we deal with in this section
 are the following: Let
$C$ be a reduced irreducible algebraic curve, and 
$$
f: C\longrightarrow C_d\tag 3.1
$$
an algebraic morphism of degree $n$
onto another reduced irreducible
algebraic curve $C_d$. Here 
we assume that everything is defined over 
the field $\C$.
We choose, once and for all,
a smooth  point $p\in C_d$ and an 
isomorphism
$$
\widehat{C_{d,p}} \overset\sim\to
\longrarrow  \Spec\big(\cring\big)
\tag 3.2
$$
of the formal completion $\widehat{C_{d,p}}$  of
$C_d$ at the divisor $p$ and the formal scheme
$\Spec\big(\C[[z]]\big)$. The isomorphism (3.2) amounts
to give a formal coordinate $z$ of $C_d$ such that the equation
$z=0$ defines the point $p$.

Let $\Cal L$ be a torsion-free sheaf on $C$ of rank 1,
and let
$$
\Cal F = f_*\Cal L
$$
be its push-forward, which is a torsion-free sheaf
of rank $n$ on $C_d$. We choose a formal trivialization
$$
\psi:\Cal F|_{\widehat{C_{d,p}}} \overset\sim\to
\longrarrow \Cal O_{\widehat{C_{d,p}}}(-p)^{\dsum n}
\tag 3.3
$$
of the vector bundle $\Cal F|_{\widehat{C_{d,p}}}$ on the
formal completion. 
We also use the same notation $\psi$ for
$$
\psi: H^0(\widehat{C_{d,p}},\Cal F|_{\widehat{C_{d,p}}} )
\overset\sim\to
\longrarrow H^0(\widehat{C_{d,p}},
\Cal O_{\widehat{C_{d,p}}}(-p)^{\dsum n})
= \big({\cring}\cdot z\big)
^{\dsum n}\;.
$$
Because of the definition, $\Cal F
= f_*\Cal L$ is an $f_*\Cal O_C$-module of
rank 1. Thus we have an $\Cal O_{C_d}$-algebra
 homomorphism
$$
f_*\Cal O_C\longrarrow \Cal End_{\Cal O_{C_d}}(\Cal F)\;.
\tag 3.4
$$
Since $\Cal L$ is torsion-free over $\Cal O_C$, 
the homomorphism (3.4) is injective.
The formal trivialization $\psi$ gives the
matrix representation of 
$$
\End(\Cal F) = H^0(C_{d},\Cal End_{\Cal O_{C_d}}(\Cal F))
\subset 
H^0(\widehat{C_{d,p}},
\Cal End_{\Cal O_{\widehat{C_{d,p}}}}(\Cal F|_{\widehat{C_{d,p}}}) )
=\End(\Cal F|_{\widehat{C_{d,p}}})
$$
 around
the smooth point $p$:
$$
\psi^*\tensor \psi: \End(\Cal F|_{\widehat{C_{d,p}}})
\overset\sim\to\longrarrow gl\big(n,{\cring}\big)\;.
$$
 In particular, the coordinate ring
of the curve $C$ admits an embedding 
$$
h:H^0(C\setminus f^{-1}(p),\Cal O_C)
\isom
H^0(C_d\setminus \{p\},f_*\Cal O_C)
\hookrarrow
H^0(C_{d}\setminus \{p\},\Cal End_{\Cal O_{C_d}}(\Cal F))
\hookrarrow gl\big(n,{\cfield}\big)
$$
into a formal loop algebra. We denote by $A$ its image:
$$
A = h\big( H^0(C\setminus f^{-1}(p),\Cal O_C)\big)
\subset gl\big(n,{\cfield}\big)\;.
\tag 3.5
$$
Let 
$$
W = \psi\big(H^0(C_d\setminus \{p\},\Cal F)\big)
\subset {\cfield}^{\dsum n}\;,
\tag 3.6
$$
and 
$$
\gamma_W:W\longrarrow 
\frac{{\cfield}^{\dsum n}}{\big({\cring}\cdot z\big)
^{\dsum n}} \tag 3.7
$$
be the natural projection. Since $\psi$ gives an 
isomorphism 
$$
H^0(\widehat{C_{d,p}},\Cal F|_{\widehat{C_{d,p}}} )
\isom \big({\cring}\cdot z\big)
^{\dsum n}\;,
$$
we have a canonical isomorphism 
$$
H^0(C_d,\Cal F)\isom \Ker(\gamma_W)\;.
\tag 3.8
$$
Similarly, the covering cohomology computation
$$
H^1(C_d,\Cal F)\isom
\frac{H^0(\widehat{C_{d,p}}\setminus\{p\},
\Cal F|_{\widehat{C_{d,p}}})}{H^0(C_d\setminus\{p\},
\Cal F)
+ H^0(\widehat{C_{d,p}},\Cal F|_{\widehat{C_{d,p}}} )}
$$
gives
$$
H^1(C_d,\Cal F)\isom \Coker(\gamma_W)\;.
\tag 3.9
$$
For more detail of this computation, see \cite{13}.
Motivated by (3.8) and (3.9), we define the 
{\it infinite-dimensional Grassmannian\/} by
$$
 Gr_n = \{W\subset {\cfield}^{\dsum n}\;|\;
\gamma_W {\text{ of (3.7) is Fredholm}}\}\;.
$$
We note that the Fredholm condition automatically
implies that $W$ is a closed subset of 
${\cfield}^{\dsum n}$
with respect to the product topology of the Krull topology
in ${\cfield}$.
The big-cell $ Gr_n ^+$ is defined to be the subset
of $ Gr_n$ consisting of points  $W$ such that 
$\gamma_W$ is an isomorphism. 

Let us denote by  $\Cal G$ the collection of 
the geometric
data 
$$
\big\langle
f:C\rarrow C_d,\Cal L, p, z, \psi
\big\rangle\tag 3.10
$$
described above. We identify 
$$
\big\langle
f':C'\rarrow C'_d,\Cal L', p', z', \psi'
\big\rangle\tag 3.11
$$
with (3.10) if there is an isomorphism
$$
j:C_d\overset\sim\to\longrarrow C'_d
$$
such that (3.10) is the pull back of (3.11) via $j$.
We have shown that a collection of
geometric data (3.10) gives rise to a collection
$$
\big\langle
i:A_d\hookrarrow A,W
\big\rangle\tag 3.12
$$
of algebraic data, where $A$ and $W$ are as in
(3.5) and (3.6), and 
$$
A_d = h(f^*(H^0(C_d\setminus \{p\},\Cal O_{C_d})))
\subset A \subset gl\big(n,{\cfield}\big)\;.\tag 3.13
$$
Of course by definition $A_d\subset\cfield$, and the 
inclusion $A_d\subset A$ is the scalar diagonal embedding.
The set of all algebraic data (3.12) is denoted by 
$\Cal A$.
We call the map
$$
\mu:\Cal G\owns \big\langle
f:C\rarrow C_d,\Cal L, p, z, \psi
\big\rangle
\longmapsto
\big\langle
i:A_d\hookrarrow A,W
\big\rangle\in \Cal A \tag 3.13
$$
the {\it cohomology map\/}, because of the
association (3.5) and (3.6).
 
{From} (3.8) and (3.9) we see that the image $\big\langle
i:A_d\hookrarrow A,W
\big\rangle$ of $\mu$
has a point $W$ of the big-cell Grassmannian
$ Gr_n ^+$ if and only if
$$
H^0(C,\Cal L) = H^1(C,\Cal L)=0\;,\tag 3.14
$$
because 
$$
H^i(C,\Cal L) \isom H^i(C_d, \Cal F)\;.
$$
Let us denote by $\Cal G^+$ the set of geometric data
satisfying (3.14), and by $\Cal A^+$ the set of algebraic
data such that $W\in Gr_n ^+$. Then we have a map
$$
\mu:\Cal G^+\longrarrow \Cal A^+\;.
$$ 
If $C$ is non-singular, then a general line bundle
$\Cal L$ on $C$ of degree $g(C)-1$ satisfies (3.14),
where $g(C)$ denotes the genus of $C$. If $C$ is 
singular, then let 
$$
r:\widetilde C\longrarrow C
$$
be the normalization of $C$. Choose a line bundle
$\widetilde{\Cal L}$ on $\widetilde C$
with vanishing cohomology groups,
and define 
$$
\Cal L = r_*(\widetilde{\Cal L})\;.
$$
Then $\Cal L$ is a torsion-free sheaf of rank 1 on
$C$ satisfying (3.14).

\bigskip\bigskip
\subhead
4. The Grassmannian and the pseudodifferential operators
\endsubhead
\medskip
In order to embed the coordinate ring $A$
of an algebraic
curve $C$ into the ring of matrix ordinary differential 
operators, we need the theory 
of matrix pseudodifferential operators. 
Let us denote by 
$$
E = \big(\C[[x]]\big)((\partial^{-1}))
$$
the set of all 
\hyphenation{pseu-do-dif-fer-en-tial}
pseudodifferential operators with coefficients in 
$\C[[x]]$, where $\partial = d/dx$. This is an 
associative algebra
and has a natural filtration 
$$
E^{(m)} = \big(\C[[x]]\big)[[\partial^{-1}]]\cdot 
\partial^m
$$ 
by the {\it order\/} of operators.
We can identify ${\cfield}$ 
with the set of pseudodifferential operators with 
constant coefficients, where
$z$ $=$ $\partial^{-1}$:
$$
{\cfield} = \C((\partial^{-1}))\subset E\;.
$$
There is also a canonical projection
$$
\rho: E \longrarrow E/Ex \isom \C((\partial^{-1})) 
= {\cfield}\;,\tag 4.1
$$
where $Ex$ is the left-maximal ideal of $E$ 
generated by $x$.
In an explicit form, this projection is given by 
$$
\rho : E\owns P = \sum_{m\in \Bbb Z} 
\partial^m\cdot a_m(x)
\longmapsto \sum_{m\in\Bbb Z} a_m(0)z^{-m}
\in {\cfield}\;.\tag 4.2
$$
It is obvious from (4.1) that ${\cfield}$ is a left 
$E$-module. The action 
is given by 
$$
P\cdot v = P\cdot\rho(Q) = \rho(PQ)\;,
$$
where
$v\in {\cfield} = E/Ex$ and $Q\in E$ is a representative 
of the equivalence
class such that $\rho(Q) = v$. 
We also use the notations
$$
\cases
D = \big(\C[[x]]\big)[\partial\/]\\
E^{(-1)} = \big(\C[[x]]\big)[[\partial^{-1}]]
\cdot\partial^{-1}\;,
\endcases
$$
which are the set of linear ordinary differential 
operators
and the set of pseudodifferential operators of 
negative order with scalar coefficients,
respectively. Note that there is a natural 
right $\big(\C[[x]]\big)$-module 
direct sum decomposition
$$
E = D\dsum E^{(-1)}\;.\tag 4.3
$$
According to this decomposition, we write $P = P^+ 
\dsum P^-$, where
$P \in E$, $P^+ \in D$, and $P^-\in E^{(-1)}$.

Now consider the matrix algebra $gl(n,E)$, 
which is the algebra of 
pseudodifferential
operators with coefficients in matrix valued functions. 
This algebra acts on the vector space
$V = {\cfield}^{\dsum n} \isom 
\big(E/Ex\big)^{\dsum n}$ from the left.
It therefore induces an infinitesimal
action of $gl(n,E)$  on the Grassmannian $Gr_n$.
The decomposition (4.3) induces 
$$
V = \C[z^{-1}]^{\dsum n}\dsum 
\big(\C[[z]]\cdot z\big)^{\dsum n}
$$
after identifying $z = \partial^{-1}$, 
and the base point $\C[z^{-1}]^{\dsum n}$ of the 
Grassmannian
$Gr_n$  is the residue class of $D^{\dsum n}$
via the projection $E^{\dsum n}
\longrarrow E^{\dsum n}\big/\big(Ex)^{\dsum n}$.
 Therefore,
the $gl(n,D)$-action on $V$ preserves $\C[z^{-1}]^{\dsum n}$.
The two theorems we need are the following:

\medskip
\proclaim{Theorem 4.1}
A pseudodifferential operator $P\in gl(n,E)$ with 
matrix coefficients
is a differential operator, i.e.~$P\in gl(n,D)$, 
if and only if
$$
P\cdot \C[z^{-1}]^{\dsum n}
\subset \C[z^{-1}]^{\dsum n}\;.
$$
\endproclaim

\medskip
\proclaim{Theorem 4.2}
Let $S\in gl(n,E)$ be a monic zero-th order 
pseudodifferential
operator of the form
$$
S = I_n + \sum_{m = 1}^\infty s_m(x) 
\partial^{-m} \;,\tag 4.4
$$
where $s_m(x)\in gl\big(n,\C[[x]]\big)$. Then the map
$$
\sigma : \Sigma\owns S\longmapsto 
W = S^{-1}\cdot \C[z^{-1}]^{\dsum n}
\in Gr_n ^+
$$
gives a bijective correspondence between the set 
$\Sigma$ of
pseudodifferential operators of the form  {\rm{(4.4)}}
and the big-cell $Gr_n ^+$ of  the Grassmannian.
\endproclaim
\medskip\noindent
Proofs of these theorems are given in 
\cite{9}.

\bigskip\bigskip
\subhead
5. Constructing commuting differential operators
from the geometric data
\endsubhead
\medskip

Once the theory of pseudodifferential operators is 
established, the passage from the geometry of curves
 to the 
ring of differential operators is straightforward.
Let
$$
\big\langle
f:C\rarrow C_d,\Cal L, p,z,\psi
\big\rangle \in \Cal G^+
$$
be a collection of geometric data with (3.14),
and let
$$
\big\langle
i:A_d\hookrarrow A,W
\big\rangle
= \mu\big(
\big\langle
f:C\rarrow C_d,\Cal L, p,z,\psi
\big\rangle\big)
$$
be the corresponding algebraic data. Then we have
$$
W\in Gr_n ^+\;,
$$
thus it corresponds to an invertible matrix
pseudodifferential operator $S\in \Sigma$
as in 
Theorem~4.2:
$$
W = S^{-1}\cdot \C[z^{-1}]^{\dsum n}\;.
\tag 5.1
$$
 Since $W$ is an $A$-module, we have
$$
A\cdot W\subset W\;. \tag 5.2
$$
The identification $z=\partial^{-1}$ makes the coordinate
ring $A$ a subalgebra of pseudodifferential operators 
with constant coefficients:
$$
A\subset gl\big(n,{\cfield}\big) = 
gl\big(n,\C((\partial^{-1}))\big)\subset gl(n,E)\;.
$$
Let
$$
B=S\cdot A \cdot S^{-1} \subset gl(n,E)
$$
be the conjugate algebra of $A$ in $gl(n,E)$ by
the operator $S$.
Then from (5.1) and (5.2), we have
$$
\aligned
B\cdot \C[z^{-1}]^{\dsum n} &=
S\cdot A \cdot S^{-1}\cdot  \C[z^{-1}]^{\dsum n}\\
&=S\cdot A\cdot W\\
&\subset S\cdot W\\
&= \C[z^{-1}]^{\dsum n}\;.
\endaligned
\tag 5.3
$$
Therefore, from Theorem~4.1, the ring $B$ consists
of differential operators:
$$
B\subset gl(n,D)\;.
$$
This completes the proof of Theorem~1.1. We remark that
the curve $C$ does not have to be irreducible for this
construction.

\bigskip\bigskip
\subhead
6. Spectral curves as an example
\endsubhead
\medskip
The theory of spectral curves by Hitchin \cite{6}
and Beauville, Narasimhan and Ramanan 
\cite{2}
provides a  method of giving a collection
of geometric data that satisfies all the conditions
of Theorem~1.1. We start with a smooth
algebraic curve $C_d$ and a  point $p$
of $C_d$. Choose a vector bundle $\Cal F$ on $C_d$
of rank $n$ and degree $n\cdot (g(C_d)-1)$
such that 
$$
H^0(C_d,\Cal F) = H^1(C_d, \Cal F) = 0\;.
$$
We also need a formal trivialization $\psi$ of $\Cal F$
as in (3.3).

Following Markman \cite{11}, let us choose a line
bundle $\Cal N$ on $C_d$ of degree at least $2g(C_d)-2$. 
The original choice is the canonical line bundle,
but any sufficiently positive $\Cal N$ will do 
for our purpose. Finally, we have to fix a {Higgs
field\/}
$$
\phi\in H^0(C_d, \Cal N\tensor \Cal End(\Cal F))\;.
$$
A Higgs field can be identified with 
a sheaf homomorphism
 $\phi : \Cal F\longrightarrow \Cal F\otimes \Cal N$, which
 induces a homomorphism 
$$ \wedge ^i\phi : \bigwedge ^i\Cal F\longrightarrow 
\bigwedge ^i\Cal F\otimes \Cal N^i\;.
$$
We call $\trace(\wedge ^i\phi )\in H^0(C_d, L^i)$
the $i$-th {\it characteristic coefficient\/} of $\phi$.
  We use the notation
 $$
ch(\phi ) = (\trace(\wedge \phi ), \trace(\wedge ^2\phi ), \cdots , 
\trace(\wedge ^n\phi ))\in 
\bigoplus _{i=1}^{n}H^0(C_d,  \Cal N ^i) \tag 6.1
$$
to denote the  characteristic coefficients
 of $\phi$.

     Consider the sheaf of symmetric algebras $Sym(\Cal N^{-1})$
 on the curve $C_d$
 generated by $\Cal N ^{-1}$. Then 
 $\Spec( Sym(\Cal N ^{-1}))$ is  the total space $|\Cal N |$ of
 the line bundle  $\Cal N$. The $i$-th characteristic coefficient 
$s_i\in H^i(C_d,\Cal N  ^i)$ gives a homomorphism
$s_i:\Cal N  ^m\longrarrow \Cal N  ^{m+i}$. 
Let $\Cal I_s$ be the sheaf of 
 ideals of $Sym(\Cal N ^{-1})$ generated by 
$$
 \bigg(1-s_1+s_2-s_3+\cdots +(-1)^n s_n\bigg)
\tensor \Cal N ^{-n} \;. \tag 6.2
$$
  Then the quotient  $Sym(\Cal N ^{-1})/\Cal I_s$ 
is a sheaf of algebras of relative
 dimension zero. We define the {\it spectral curve\/} by
 $$
C_s =\Spec\big( Sym(\Cal N ^{-1})/\Cal I_s\big)\;. \tag 6.3
$$
In a more geometric language, the spectral curve $C_s$ is a subvariety
of $|\Cal N |$,   where
$|\Cal N |$ is considered to be an open algebraic surface.
The projection of $|\Cal N |$ onto $C_d$ defines a natural map 
$$
f_s:C_s\longrarrow C_d
$$ 
of degree $n$. Its fiber at $q\in C_d$
is the set of solutions of the polynomial equation 
$$
\bigg(1-s_1+s_2-s_3+\cdots +(-1)^n s_n\bigg)\tensor y^n = 
y^n -\bar s_1 y^{n-1}+\bar s_2
y^{n-2} - \cdots +(-1)^n \bar s_n =0
$$
evaluated at $q$,
where $y\in \Cal N  ^{-1}$ is a linear coordinate of the fiber of
$\Cal N$, and 
$$
\bar s_i = s_i \tensor y^i \in \Cal N  ^i \tensor
\Cal N  ^{-i} = \Cal O_{C_d}
$$ 
is a regular function defined locally near $q\in C_d$.
This polynomial
equation is also the characteristic equation of the homomorphism
$\phi : \Cal N  ^{-1} \longrarrow \Cal End(\Cal F)$ 
at $q$ when $ch(\phi) = s$. 
 Note that the homomorphism 
$$
\phi:\Cal N^{-1}\longrightarrow \Cal End(\Cal F)
$$
 induces an algebra homomorphism 
$$
Sym(\Cal N^{-1})\longrightarrow
 \Cal End(\Cal F)\;.
$$
 This algebra homomorphism factors
 through 
$\Cal O_{C_s}=Sym(\Cal N^{-1})/\Cal I_s$ if $ch(\phi )=s$, which
  gives a $Sym(\Cal N^{-1})/\Cal I_s$-module
structure in $\Cal F$. Since the rank
of $Sym(\Cal N^{-1})/\Cal I_s$ over $\Cal O_{C_d}$ is $n$, 
$\Cal F$ defines a line bundle over $C_s$, 
which we denote by $\Cal L_s$.  It is
 clear that $({f_s})_*(\Cal L_s)=\Cal F$.
Thus we have constructed a desired geometric data
$$
\big\langle
f_s:C_s\rarrow C_d,\Cal L_s,p,z,\psi
\big\rangle\in \Cal G^+\;.
$$

\bigskip\bigskip
\subhead
7. {From} differential operators to geometry
\endsubhead
\medskip
Let $B\subset gl(n,D)$ be an integral domain that
satisfies the condition of Theorem~1.2. It is
easy to show that
\medskip
\proclaim{Lemma 7.1}
There is an element $S\in \Sigma$
such that
$$
P = S\cdot I_n\cdot \partial^r \cdot S^{-1}\;,
$$ 
where $r$ is the order of $P$ as a monic elliptic ordinary 
differential operator. 
\endproclaim
\medskip\noindent
The 
operator $S$ relates the ring of differential
operators to the formal loop algebras: 
$$
A = S^{-1}\cdot B\cdot S\subset gl\big(n,{\cfield}\big)\;,
$$
where $z=\partial^{-1}$. This is
 because every element of $A$
commutes with $I_n\cdot \partial^{r}$, and hence it has 
constant coefficients. Similarly, we define
$$
A_d = S^{-1}\cdot B_d\cdot S\subset gl\big(n,{\cfield}\big)\;.
$$
Since $B_d\subset \C((P^{-1/r}))$, we have
$$
A_d\subset I_n\cdot {\cfield}\;,
$$
i.e., it is a scalar diagonal subalgebra of $A$.
By definition, the pair $(A,A_d)$ satisfies 
Condition~2.1. Thus by Theorem~2.2, we have
a degree $n$ morphism
$$
f:C\longrarrow C_d\;,
$$
where $C_d$ is a one-point completion of
$\Spec(A_d)$ with a smooth point $p$. 

We still have to construct a torsion-free
rank one sheaf $\Cal L$ on $C$. Note that
$D^{\dsum n}$ is a left $gl(n,D)$ right $\C[[x]]$
bimodule. We consider $D^{\dsum n}$ as bimodule 
over $B\tensor \C[[x]]$, accordingly. 
Similarly, $I_n\cdot D$ is a bimodule over
$B_d\tensor \C[[x]]$. It is known \cite{15}
that this module structure 
gives rise to a rank one sheaf on $C_d\times \Spec(\C[[x]])$.
Since $D^{\dsum n}$ is rank $n$ over $D$ and $B$ is rank $n$ over
$B_d$, we conclude that $B\tensor \C[[x]]$-bimodule
$D^{\dsum n}$ defines a torsion-free rank one sheaf on 
$C\times \Spec(\C[[x]])$. Restricting it to
$C = C\times \{0\}$, we obtain 
a torsion-free rank one sheaf $\Cal L$ defined
on $C$. The vanishing  of the cohomology groups (3.14) can be
proved in a similar method developed in \cite{9} 
and\cite{13}. This completes the proof of Theorem~1.2.

\bigskip\bigskip
\def\nspace{\lineskip=1pt\baselineskip=12pt\lineskiplimit=0pt}
\def\references#1#2{\bigskip\par\centerline{\bf #1}\medskip
                    \parindent=#2 pt\nspace}
\def\Ref[#1]{\par\smallskip\hang\indent\llap{\hbox to\parindent
              {[#1]\hfil\enspace}}\ignorespaces}
\references{Bibliography}{20}
\medskip\noindent
\Ref[1] M.~R.~Adams and M.~J.~Bergvelt: 
The Krichever map,
vector bundles over algebraic curves, 
and Heisenberg algebras,
Commun.\ Math.\ Phys.\ {\bf 154} (1993)  265--305.
\Ref[2] A.~Beauville, M.~S.~Narasimhan and 
S.~Ramanan: Spectral
curves and the generalized theta divisor, 
Journ.\ Reine Angew.\ Math.\
{\bf 398} (1989) 169--179.
\Ref[3] J. L.~Burchnall and T. W.~Chaundy: Commutative ordinary
differential operators, Proc.\ London Math.\ Soc. Ser.~2, {\bf 21}
(1923) 420-440; Proc.\ Royal Soc.\ London Ser.~A, {\bf 118} (1928)
557--583.
\Ref[4] R.~Donagi and E.~Markman: Spectral curves, 
algebraically completely integrable Hamiltonian systems, and
moduli of bundles, alg-geom/9507017
(http://xxx.lanl.gov/ list/alg-geom/9507)
\Ref[5] I.~M.~Gel'fand and L.~A.~Dickey: Asymptotic behavior of the
resolvent of Sturm-Liouville equations and the algebra of the
Korteweg-de~Vries equations, Russ.\ Math.\ Surveys, {\bf 30}
(1975) 77--113; Fractional powers of operators and Hamiltonian
systems, Func.\ Anal.\ Appl.\ {\bf 10} (1976) 259--273.
\Ref[6] N.~Hitchin: Stable bundles and integrable systems,
Duke Math.\ J.\ {\bf 54} (1987) 91--114.
\Ref[7] M.~Kimura: Ph.\ D.\ Thesis, University of California,
Davis (1996)
\Ref[8] I.~M.~Krichever: Methods of algebraic geometry in the theory of
nonlinear equations, Russ.\ Math.\ Surveys {\bf 32} (1977) 185--214.
\Ref[9] Y.~Li and M.~Mulase: 
Category of morphisms of 
algebraic curves and a characterization of Prym varieties,
alg-geom/9203002 (http://xxx.lanl.gov/list/alg-geom/9203).
\Ref[10] Y.~Li and M.~Mulase: Hitchin
systems and KP equations, to appear 
\Ref[11] E.~Markman: 
Spectral curves and integrable systems,
 Compositio Math.\ {\bf 93} (1994) 255--290. 
\Ref[12] M.~Mulase: Cohomological 
structure in soliton
equations and Jacobian varieties, 
J.\ Differ.\ Geom.\ {\bf 19} (1984)
403--430.
\Ref[13] M.~Mulase: Category of vector bundles
on algebraic curves and infinite-dimensional 
Grassmannians, 
Intern.\ J.\ of Math.\ {\bf 1} (1990)   293--342.
\Ref[14] M.~Mulase:
Algebraic theory of the KP equations,
in {\it Perspectives in Mathematical Physics\/}, 
R.~Penner and S.~T.~Yau, Eds., Intern.\ Press Co., (1994) 157--223.  
\Ref[15] D.~Mumford: An algebro-geometric constructions of commuting
operators and of solutions to the Toda lattice equations,
Korteweg-de~Vries equations and related non-linear equations, 
in Proc.\
Internat.\ Symp.\ on Alg.\ Geom., Kyoto 1977, Kinokuniya Publ.\ (1978)
115--153.
\Ref[16] E.~Previato and G.~Wilson: Vector bundles 
over curves and solutions of the KP equations, 
Proc.\ Symposia in Pure Math.\ {\bf 49} (1989) 553--569.
\Ref[17] M.~Sato: Soliton equations as dynamical systems on an infinite
dimensional Grassmannian manifold, Kokyuroku, Res.\ Inst.\ Math.\ Sci.,
Kyoto Univ. {\bf 439} (1981) 30--46.
\Ref[18] I.~Schur: \"Uber vertauschbare lineare
Differentialausdr\"ucke,  Sitzungsber.\ der Berliner Math.\
Gesel.\ {\bf 4} (1905) 2--8.
\Ref[19] G.~B.~Segal and G.~Wilson: Loop groups and equations of $KdV$
type, Publ.\ Math.\ I.H.E.S. {\bf 61} (1985) 5--65.
\Ref[20] J.-L.~Verdier: Equations differentielles alg\'ebriques,
S\'eminaire de l\'Ecole Normale \newline Sup\'erieure 1979-82, 
Birkh\"auser (1983) 215--236.
\Ref[21] G.~Wallenberg: \"Uber die Vertauschbarkeit homogener
linearer  Differentialausdr\"ucke, Archiv der Math.\ u.\ Phys.,
Drittle Reihe {\bf 4} (1903) 252--268.
\vfill
\eject
\enddocument
\end